\documentclass{amsart}
\usepackage{amscd}
\usepackage{a4wide}
\usepackage{amssymb,amsmath}  
\usepackage{latexsym} 
\newcommand{\Ps}{\mathbf{P}}

\newcommand{\Z}{\mathbf{Z}}
\newcommand{\ra}{\rightarrow}
    \newtheorem{Lem}{Lemma}
    \newtheorem{Prop}[Lem]{Proposition}
    \newtheorem{Thm}[Lem]{Theorem}
    \newtheorem{Cor}[Lem]{Corollary}

   \theoremstyle{definition}
    \newtheorem{Def}[Lem]{Definition}
    
    \newtheorem{Que}[Lem]{Question}
    
     \theoremstyle{remark}
    
    \newtheorem{Rem}[Lem]{Remark}
    
  \DeclareMathOperator{\mult}{mult}

 \DeclareMathOperator{\cha}{char}

 \DeclareMathOperator{\Aut}{Aut}

\begin{document}
\title[Geometry of the Weierstrass scheme]{Locally trivial families of hyperelliptic curves: the geometry of the Weierstrass scheme}
\author{Remke Kloosterman}
\email{r.n.kloosterman@math.rug.nl}
\address{Department of Mathematics and Computer Science\\ University of Groningen\\ PO Box 800\\ 9700 AV  Groningen\\ The Netherlands}

\author{Orsola Tommasi}
\email{tommasi@math.ru.nl}
\address{Department of Mathematics\\ Radboud University Nij\-me\-gen\\ Toernooiveld 1 \\ 6525 ED Nij\-me\-gen \\ The Netherlands}
\begin{abstract}
In this paper we describe some geometrical properties of the Weierstrass scheme of locally trivial hyperelliptic fibrations.
\end{abstract}
\maketitle

\section{Introduction}
In this introduction suppose $K$ is a field of characteristic 0.

Starting point of the research on which we report in this paper was a discussion between the second author and Emilia Mezzetti, concerning the solution of Mezzetti and Portelli to the following problem:

\begin{Que}Fix a point $p\in\Ps^2$ and a reduced cubic curve $C\subset \Ps^2$ such that $p\not \in C$. Let $L$ be a line through $p$. Suppose $C\cap L$ consists of three distinct points, say $p_{L,1},p_{L,2},p_{L,3}$. 

For which pairs $(C,p)$ does an enumeration of the elements in $C\cap L=\{p_{L,1},p_{L,2},p_{L,3}\}$ exist, such that for almost all lines $L$ through $p$ the cross-ratio of $p,p_{L,1}, p_{L,2}$ and $p_{L,3}$ is independent of $L$? \end{Que} 

The answer to this question is: 
\begin{Thm}\label{thm2} $(C,p)$ has the above mentioned property if and only if one of the following holds: 
\begin{enumerate}
\item  $C$ is the union of three lines passing through one point and $p$ is a point outside $C$.
\item  $C$ is the union of a non-singular conic $Q$ and a line $L$, not tangent to the conic, and $p$ is the polar point of  $L$ with respect to $Q$.
\item $C$ is a cuspidal cubic, and $p$ is the intersection point of the flex line of $C$ and the tangent line at the cusp.
\item $C$ is a smooth elliptic curve with $j$-invariant 0, and $p$ is the intersection point of the tangent lines at three collinear flexes.
\end{enumerate}
\end{Thm}

One can prove this in several ways. 
One strategy, using the theory of elliptic surfaces, is explained in Remark~\ref{remell}. 
We give also another proof of this Theorem, which applies to a generalized setting.
Indeed, there is no reason to restrict to cubic curves, or to assume $p\not \in C$. We consider pairs $(C,p)$ with $\deg(C)-\mult_p(C)=d\geq 3$, such that $p$ is not on a line contained in $C$,  with the property that for every two lines $L$ and $M$ through $p$ such that $\# (C-\{p\}) \cap L = \# (C-\{p\}) \cap M=d$, we can write $(C-\{p\}) \cap L=\{p_1,\dots,p_d\} $ and $(C-\{p\}) \cap M=\{q_1,\dots,q_d\}$ in such a way that the 
lines $p_iq_i$ are concurrent. If $(C,p)$ has this property then we say that $(C,p)$ has \emph{constant moduli}.

The main result is the following:
\begin{Thm} \label{mainthm} Fix $p\in \Ps^2$. Suppose $C$ is a reduced plane curve of degree $d+m$, with $m$ the multiplicity of $C$ at $p$. Suppose no line through p is contained in $C$. 
Then the following are equivalent: 
\begin{enumerate}
\item The pair $(C,p)$ has constant moduli.
\item There exists an automorphism (if $m=0$) or a birational automorphism (if $m>0$) $\varphi$ of $\Ps^2$ and a homogenous polynomial $H\in K[Y,Z]$ of degree $n$, such that $\varphi(p)=[1:0:0]$, almost every line through $p$ is mapped to a line through $[1:0:0]$ and the closure of $\varphi(C)$ is the zero-set of
\[ \prod^{d/k}_{t=1} (X^k Z^{n-k}-\alpha_t H(Y,Z)) \mbox{ or } X \prod_{t=1}^{(d-1)/k} (X^kZ^{n-k}-\alpha_t H(Y,Z)), \]
where $k$ is an integer dividing either $d$ or $d-1$, and $\alpha_t \in \overline{K}$ for every index $t$.
\end{enumerate}
\end{Thm}

The equations for $(C,p)$ give the following:
\begin{Cor} With notation as above, if $(C,p)$ has constant moduli then 
for some positive integer $k$, 
the group $G$, with
\[ G \cong\frac{\Z}{k\Z} \times \frac{\Z}{\frac{d}{k}\Z} \mbox{ or } G \cong  \frac{\Z}{k\Z} \times \frac{\Z}{\frac{d-1}{k}\Z} \]
is a subgroup of $\Aut(C)$. \end{Cor} 

To every pair $(C,p)$ as in the hypotheses of Theorem~\ref{mainthm} we can associate a surface fibered in hyperelliptic curves of genus either $(d-1)/2$ or $(d-2)/2$. This will be explained explicitly in Remark~\ref{remell}, in the case of curves of degree $3$.
Basically, we take for almost every line $L$ through $p$ the double cover of $L$ ramified at $L\cap C$ or at $L\cap (C-\{p\})$, depending on whether $d$ is even or odd. The pair $(C,p)$ has constant moduli if and only if the associated family of hyperelliptic curves is locally trivial.

In fact, it is true (see Proposition~\ref{surf}) that every surface with a hyperelliptic fibration over $\Ps^1$, admitting either a two-section invariant under the hyperelliptic involution or a section, is birational to the surface associated to some pair $(C,p)$.
As a consequence, Theorem~\ref{mainthm} provides a local form for locally trivial families of hyperelliptic curves. 

\section{Proof of the theorem}
In this section we prove Theorem~\ref{mainthm}. 

We always follow the convention that the curve $C$ has degree $d+m$, where $m\geq 0$ is the multiplicity of $C$ at $p$.
Moreover, we assume that $d>2$ and $\cha(K)>d$ or $\cha(K)=0$. 

\begin{Def} Let $L$ and $M$ be distinct lines through $p$;
 let $A\subset L$ and $B\subset M$ be two sets of $k$ points. These two sets have the same \emph{moduli} if there is an isomorphism $L\longrightarrow M$ mapping $A$ to $B$ and fixing $L\cap M$. 

We say that the pair $(C,p)$ has \emph{constant moduli} if for almost all lines $L$ and $M$ through $p$ the sets $(C-\{p\}) \cap M$ and $(C-\{p\}) \cap L$ have the same moduli.
\end{Def}

Two sets of collinear points $A$, $B$ (not on the same line) 
have the same moduli if and only if one can write $A=\{p_1,\dots, p_k\}$ and $B=\{q_1,\dots,q_k\}$ in such a way that all 
lines $p_iq_i$ (with $p_i\neq q_i$) have a common intersection point. 
The reader can easily verify this, e.g., by choosing coordinates on the lines containing $A$, respectively, $B$.

\begin{Lem}\label{remtang} Suppose $(C,p)$ has constant moduli. Suppose $L$ is a line such that $S=(C-\{p\}) \cap L$ consists of $d$ points. Then the $d$ tangent lines of $C$ at the points in $S$ are concurrent.\end{Lem}

\begin{proof}Write $S=\{p_1,\dots,p_d\}$. For almost every line $M$ through $p$, we can write $(C-\{p\})\cap M=\{q_1,\dots,q_d\}$ in such a way that 
the $p_iq_i$ are concurrent.   Moving $M$ to $L$ shows that the $d$ tangent lines pass through a unique point.
\end{proof}

\begin{Lem} 
Suppose that $(C,p)$ has constant moduli.

Suppose $\ell$ is a line through $p$ such that $n=\# (C-\{p\}) \cap \ell <d$.
Then $n=0$ or $n=1$.
\end{Lem}

\begin{proof} Take a line $\mu$ through $p$ intersecting $C-\{p\}$ in exactly $d$ points. By a limit position argument as in the proof of Lemma~\ref{remtang}, there exist lines $\ell_1,\dots,\ell_d$ connecting all the points in $(C-\{p\})\cap\mu$ with the points in $(C-\{p\})\cap\ell$, such that $\ell_1,\dots,\ell_d$ have a common intersection point. Since $(C-\{p\})\cap\ell$ has less than $d$ points, the common intersection point lies on $\ell$, hence $(C-\{p\})\cap\ell$ is empty or consists of one point.\end{proof}

We describe next some reduction steps. 

\begin{Rem} \label{RemRed} Let $C$ be a reduced plane curve. Fix some $p\in 
\Ps^2$, not on a line contained in $C$.

Without loss of generality we may assume that $p=[1:0:0]$. If we set $Z=1$ 
then an affine equation of $C$ is of the form
\[ \sum_{k=0}^{d}  f_{k}(y)x^k =0. \]
After multiplying this equation by $(f_{d})^{d-1}$ and replacing $x$ by $x/f_{d}$ we may assume that $f_{d}=1$. Then replacing $x$ by $x-\frac{1}{d} f_{d-1}$ allows us to assume that $f_{d-1}=0$.

Both substitutions define birational automorphisms of $\Ps^2$. 
The first map may not be defined on 
lines  $L:aY=bZ$ such that $(C.L)_p>\mult_p(C)$. The second map may not be 
defined on the line $Z=0$. 
Almost all lines through $p$ are mapped to lines through $p$. 
Note that it can happen that, as a result of these substitutions, $C$ acquires a line through $p$ as a component, possibly with multiplicity. In this case we choose to discard this component, because lines through $p$ are not relevant for the property of having constant moduli. If we discard after each birational map these lines through $p$ then both maps are defined on a dense open subset of $C$.

If $p \not \in C$, then both birational maps are isomorphisms. 

Note that these substitutions do not alter $d$, but can change the value of $m$. 
Also the singularities of $C$ can change. The reason why we study the above type of equation is that surfaces with an hyperelliptic fibration are usually written in the following form, which is a higher-degree analogue of the Weierstrass form:
\[ z^2=x^d+f_{d-2}(y) x^{d-2}+\dots + f_0(y). \]
\end{Rem}

In view of the preceding remark, from now on we assume $p=[1:0:0]$ and $C$ is a reduced curve of degree $d+m$, such that $C$ is the zero-set of 
\[\tag{*} G(X,Y,Z)=Z^{d-m} X^d+\sum_{k=0}^{d-2} F_{d+m-k}(Y,Z) X^k,\]
with $F_h$ a homogeneous form of degree $h$ in $Y$ and $Z$, and $Z\nmid G$.

\begin{Def}
Suppose that $(C,p)$ has constant moduli. 
A \emph{special line} of $C$ is a line $\ell$ through $p$, such that $(C-\{p\}) \cap L$ has 
one element. This point is called a \emph{special point}.

Suppose $\ell$ is a non-special line through $p$. Then we denote by $T_\ell$ the intersection point of all tangent lines to $C$ at points in $(C-\{p\})\cap\ell$.  
\end{Def}

\begin{Prop}\label{PrpTline} Suppose $(C,p)$ has constant moduli.
Then 
\[T:=\overline{ \{T_\ell \mid \ell \mbox{ line through } p \mbox{ and } \# (C-\{p\})\cap\ell =d \} } \]
is a point or a line.

\end{Prop}

\begin{proof} Fix a general line $L$ of the form $Y=y_0 Z$. From Lemma~\ref{remtang} it follows that the lines  tangent to  $C$ at the points in $(C-\{p\})\cap L$ have a common intersection point $q$. 
  Possibly after an automorphism of the form $Y\mapsto Y+cZ$, we may assume that $q=[\bar{x}:\bar{y}:1]$, hence  for every point $[x_0:y_0:1] \in C(\bar{K}) \cap L$ we have 
\[ \frac{\partial{G}}{\partial{X}} (x_0,y_0,1) (x_0-\bar{x}) +\frac{\partial{G}}{\partial{Y}} (x_0,y_0,1)(y_0-\bar{y})  =0.\]
This gives rise to a polynomial in $x_0$ of positive degree, with the same zeros as $G(x_0,y_0,1)$ considered as polynomial in $x_0$. Comparing the two highest coefficients shows that $\bar{x}=0$, so $T$ is contained in a line. Since there is a dominant map from $\Ps^1$ to $T$, we obtain that $T$ is either a single point or the whole line $X=0$.

\end{proof}

\begin{Rem}
Suppose $(C,p)$ has constant moduli. Then Theorem 3 implies that all special points are contained in $T$. 
\end{Rem}

This yields a geometrical construction for finding the center $p$ of the pencil, once the curve $C$ is known. If the curve $C$ is non-singular outside $\{p\}$, the line $T$ is a line on which $d+m$ flex points lie. Moreover, the tangent lines at those $d$ flexes intersect $C$ with multiplicity $d$, and have a common point, which is the point $p$. 

\begin{proof}[Proof of Theorem~\ref{mainthm}] (1)$\Rightarrow$(2).
 
By Remark~\ref{RemRed} we may assume that the equation of $C$ is of the form (*).
Suppose first that $C$ does not contain the line $T: X=0$.

>From the proof of Proposition~\ref{PrpTline} we obtain for almost all $y_0\in\bar K$ the relation 
\[ \frac{\partial{G}}{\partial{X}} (x,y_0,1)(x-\bar{x}) +\frac{\partial{G}}{\partial{Y}} (x,y_0,1)(y_0-\bar{y})=d G(x,y_0,1).\]
Comparing coefficients gives $\bar{x}=0$ and
\[ (y_0-\bar y)\frac{\partial F_{m+h}}{\partial Y}(y_0,1) = h F_{m+h}(y_0,1),\]
for all integers $h$ such that $2\leq h\leq d$. 
This implies that
\[h {F_{m+h}(y_0,z)}\left({\frac {\partial F_{m+h}}{\partial Y}}(y_0,z)\right)^{-1}\]
is independent of the index $h$, provided $\frac{\partial F_{m+h}}{\partial Y}(y_0,z)$ is different from 0. 

Let $n:=\gcd\{h \mid F_{m+h}(Y,Z)\neq 0\}$. From the above it follows that for all $h,j$ such that $F_{m+h},F_{m+j}\neq 0$ we have 
\[  F_{m+h}^{j}=c_{h,j} Z^{m(j-h)} F_{m+j}^{h}, \] 
hence there exists a homogenous polynomial $H(Y,Z)$  such that $F_{m+h}(y,1)=\lambda_{h/n} H^{h/n}(y,1)$ if $n$ divides $h$ and $2\leq h \leq d$ and $F_{m+h}=0$ otherwise. Thus
\[ G := \sum_{t=0}^{d/n} \lambda_{t} (H(Y,Z))^t (X^n)^{d/n-t} Z^{(d-tn)m/d}= \alpha_0 \prod_t (X^nZ^{nm/d}-\alpha_t H(Y,Z))\] 
for some $\alpha_t \in \bar{K}$.

In the case that $C$ contains $T$, consider $\overline{C-T}$. This curve has constant moduli, and hence has an equation of the above form. Multiplying this equation with the equation for $T$ gives an equation as in (2).

(2)$\Rightarrow$(1).

Suppose $C$ is of the form
\[\prod_t F_t =0,\;\;\;\;F_t(X,Y,Z):= X^nZ^{nm/d}-\alpha_t H(Y,Z). \]
Denote by $C_t$ the curve of equation $F_t=0$, and define $\sigma : [X:Y:Z] \mapsto [\zeta_n X:Y :Z ]$, $\zeta_n$ a primitive $n$-th root of unity.

Fix some $t$ and two lines $\ell_1,\ell_2$ through $[1:0:0]$ both intersecting $C_t$ in $n$ distinct points.
Let $L_n$ be a line connecting one point of $\ell_1\cap C$ and one point of $\ell_2\cap C$.
Pose $L_k=\sigma^k(L_n)$ for $k=1,\dots,n-1$. Then $L_1,\dots,L_n$ are lines connecting the points in $C\cap\ell_1$ with those in $C\cap\ell_2$. Note that $L_1\cap\dots\cap L_n$ is one point, lying on the line $X=0$.

The morphism $[X:Y:Z] \mapsto [\sqrt[n]{\alpha_{t+1}/\alpha_t}X:Y:Z]$ maps $C_t$ to $C_{t+1}$, maps $L_k$ to a line connecting a point of $C_{t+1}\cap \ell_1$ with a point of $C_{t+1}\cap \ell_2$ and fixes the common intersection point of the $L_i$.

If $C$ is of the form $X\prod_t F_t=0$, then $C$ is the union of the line $X=0$ and a curve $C'$ of equation $\prod_t F_t=0$. We already know that the pair $(C',p=[1:0:0])$ satisfies condition (1). 
Then the claim follows from the fact that the lines connecting the intersections of $C'$ with two general lines of the pencil with center $p$ intersect always at a point on $X=0$. 
\end{proof}

We are ready to prove Theorem~\ref{thm2}.

\begin{proof}[Proof of Theorem~\ref{thm2}] 
We are in the case $d=3,m=0$.
We know that $(C,p)$ has constant moduli if and only if there exists an automorphism of $\Ps^2$ mapping $p$ to $[1:0:0]$ and $C$ to the zero-set of a polynomial of the form 
\[ X^3+F_3(Y,Z), \;\;\; X(X^2+F_2(Y,Z)) \mbox{ or } X(X+F_1(Y,Z))(X+\lambda F_1(Y,Z)) \] 
with $F_i(Y,Z)$ a non-zero homogenous polynomial of degree $i$.

In all three cases, if $F_i=0$ defines one point in $\Ps^1$, then $C$ is the union of three concurrent lines.
If $F_2=0$ defines two distinct points in $\Ps^1$ then $X(X^2+F_2(Y,Z))$ is the union of a conic and a line, and $p$ is the polar point of $L$ with respect to $C$.

If $F_3=0$ defines two distinct points in $\Ps^1$ then $X^3+F_3(Y,Z)$ is a cuspidal cubic and $p$ the intersection point of the tangent line at the cusp and the flex line of $C$.

If $F_3=0$ defines three distinct points in $\Ps^1$ then $X^3+F_3(Y,Z)$ is a smooth elliptic curve with $j$-invariant 0, and $p$ is the intersection point of the tangent lines at three collinear flexes.
\end{proof}

\begin{Rem}\label{remell}Let $C$ be a cubic curve and $p$ be a point not on $C$. Let $X_1$ be the blow-up of $\Ps^2$ in $p$. Denote by $E$ the exceptional divisor.  
Fix some fiber $F$ of the ruling of $X_1$. Let $X_2$ be the surface obtained from $X_1$ by blowing it up along $E\cap F$ and contracting the strict transform of $F$. Note that $X_2$ is the second Hirzebruch surface (see \cite[V.4]{BPV}). 
Let $D$ be the union of the strict transforms of $C$ and $E$ on $X_2$. 
Consider the double cover $Y$ of $X_2$ ramified along $D$. The surface $Y$ is a rational elliptic surface. Then we have that $(C,p)$ has constant moduli if and only if the $j$-invariants of all non-singular fibers of $Y\rightarrow \Ps^1$ coincide.

In the above process we have to choose a fiber $F$. We may assume that $F$ intersects the strict transform of $C$ in three distinct points. After fixing the morphism $X_1 \ra \Ps^1$, we may assume that $F$ lies above the point $[1:0]$. 
Then $Y$ has an affine equation of the form
\[ y^2=x^3+f_1(t)x^2+f_2(t)x+f_3(t) \]
with $\deg f_i \leq i$. After substituting $x=x'-\frac{1}{3}f_1(t)$ we may assume that $f_1=0$. Since $F$ intersects $C$ in three distinct points we have $3\deg(f_2)+2\deg(f_3)=6$.

It is known (\cite[Section III.1]{BU}) that 
the fibers of the elliptic fibration on $Y$ have constant $j$-invariant if and only if $f_2^3=cf_3^2, c\in K$ or $f_2f_3=0$. 
This gives another proof of Theorem~\ref{thm2}.
\end{Rem}

The equation given in Theorem~\ref{mainthm} has the following geometric interpretations in the case $d=4, m=0$: 

\begin{Cor} Suppose $d=4,m=0$. Then $(C,p)$ has constant moduli if and only if one of the following holds: 
\begin{enumerate}
\item $C$ is the union of four concurrent lines. 
\item $C$ is the union of two conics $C_i$. Then the tangent lines of $C_1$ and $C_2$ at the intersection points have to coincide, and $p$ is the common polar point.
\item $C$ is the union of a cubic $E$ and a line $L$. Then $E$ is either cuspidal or smooth with $j(E)=0$, and $L$ is the line passing through the flex and the cusp or the line passing through three flexes.
\item $C$ is a cyclic cover of $\Ps^1$ of degree 4 ramified in 4 points, and $p$ is the intersection point of the tangents at $C$ at the four ramification points.
\item $C$ is irreducible, has two 4-flexes and a tacnode and all these points are collinear. The intersection point of the tangent line at the tacnode and two flex lines of the $4$-flexes is $p$. The normalization of $C$ is an elliptic curve with $j$-invariant 1728.

\item $C$ is irreducible and has a triple point of the form $y^3=x^4$ and a 4-flex. The point $p$ is the intersection point of the tangent line at the triple point and the 4-flex line. The normalization of $C$ has genus 0.
\end{enumerate}
\end{Cor}

We prove now the assertion of the Introduction about the equivalence of considering pairs $(C,p)$ with constant moduli and locally trivial hyperelliptic fibrations over $\Ps^1$ admitting a section or a two-section invariant under the hyperelliptic involution.

\begin{Prop}\label{surf}Suppose we have a family of hyperelliptic curves $\pi:Y\ra \Ps^1$, together with a section $\sigma:\Ps^1\ra Y$, or a two-section invariant under the hyperelliptic involution.

Then we can construct a pair $(C,p)$ such that its associated surface (as in the Introduction) is $Y$.
\end{Prop}

\begin{proof} Suppose we have a family of hyperelliptic curves $\pi:Y\ra \Ps^1$, together with a section $\sigma:\Ps^1\ra Y$, such that $\sigma(\Ps^1)\cap\pi^{-1}(t)$ is a Weierstrass point of $\pi^{-1}(t)$ for almost all $t$. Let $\tilde{Y}$ be the surface obtained from $Y$ by contracting all components of fibers not intersecting $\sigma(\Ps^1)$. Consider the quotient $X=\tilde{Y}/\langle \iota\rangle$, with $\iota$ the hyperelliptic involution on all fibers. Then $X$ is a ruled surface. Let $\tilde C$ be the image of the fixed locus of the hyperelliptic involution minus the image of $\sigma(\Ps^1)$. Since $X$ is a ruled surface we can construct a rational map $\psi: X \dashrightarrow \Ps^2$ such that the ruling induces a pencil of lines through the point $p$, with $\{p\}=\psi(\sigma(\Ps^1))$. Then we take $C$ to be the closure of $\psi(\tilde{C})$. 

If the section is not a Weierstrass point for almost all $t$, then the section and its conjugate form a two-section invariant under the hyperelliptic involution. The reasoning above can easily be adapted to the case of an invariant two-section.
\end{proof}

\begin{Rem}
Theorem~\ref{mainthm} gives an indication of 
how locally trivial families of hyperelliptic curves degenerate. If the $2g+2$ Weierstrass points $\{P_1,\dots, P_{2g+2}\}$ are in `general' position, then there is locally a unique way to degenerate this family. (In the case of elliptic curves this result is known and formulated as follows: Suppose $\pi: X \ra C$ is an elliptic surface with constant $j$-invariant different from 0,1728 then all singular fibers are of type $I_0^*$.)

If
\[ \{ \varphi \in \Aut(\Ps^1) \mid \varphi(\{P_1,\dots, P_{2g+2}\}) =\{P_1,\dots, P_{2g+2}\} \}\]
has more than one element then the position of the Weierstrass points (i.e., the moduli of the hyperelliptic curve) does not determine the degeneration locally, as we have already seen in the case of elliptic surfaces with constant $j$-invariant 0 and 1728. 

Moreover, it seems that one has to add one global combinatorial condition to describe all configurations of degenerations possible in a complete family. (In the case of elliptic curves, the sum of the Euler numbers of the singular fibers equals $12(p_g(X)+1)$.)
\end{Rem}

\end{document}